\documentclass[11pt]{amsart}

\usepackage{amsmath}
\usepackage{amssymb}
\usepackage{mathrsfs}
\usepackage[all]{xy}

\newtheorem{theorem}{Theorem}[section]

\newtheorem{lemma}[theorem]{Lemma}

\newtheorem{corollary}[theorem]{Corollary}

\theoremstyle{definition}
\newtheorem{definition}[theorem]{Definition}

\newtheorem{question}[theorem]{Question}

\theoremstyle{remark}
\newtheorem{remark}[theorem]{Remark}

\newcount\skewfactor
\def\mathunderaccent#1#2 {\let\theaccent#1\skewfactor#2
\mathpalette\putaccentunder}
\def\putaccentunder#1#2{\oalign{$#1#2$\crcr\hidewidth
\vbox to.2ex{\hbox{$#1\skew\skewfactor\theaccent{}$}\vss}\hidewidth}}
\def\name{\mathunderaccent\tilde-3 }

\def\smallbox#1{\leavevmode\thinspace\hbox{\vrule\vtop{\vbox
   {\hrule\kern1pt\hbox{\vphantom{\tt/}\thinspace{\tt#1}\thinspace}}
   \kern1pt\hrule}\vrule}\thinspace}

\newcommand{\cf}{{\rm cf}}

\def\qedref#1{$\qed_{\reforiginal{#1}}$}

\title[Configurations]{On configurations concerning cardinal characteristics at regular cardinals}
\author{Omer Ben-Neria}
\address{Einstein Institute of Mathematics
 The Hebrew University of Jerusalem,
 Jerusalem 91904, Israel}
\email{omer.bn@mail.huji.ac.il}

\author{Shimon Garti}
\address{Einstein Institute of Mathematics
 The Hebrew University of Jerusalem,
 Jerusalem 91904, Israel}
\email{shimon.garty@mail.huji.ac.il}

\thanks{The first author was partially supported by the Israel Science Foundation grant no. 1832/19. The second author was supported by ERC grant no. 338821.}

\begin{document}
\let\labeloriginal\label
\let\reforiginal\ref
\def\ref#1{\reforiginal{#1}}
\def\label#1{\labeloriginal{#1}}

\begin{abstract}
We study the consistency and consistency strength of various configurations concerning the cardinal characteristics $\mathfrak{s}_\theta, \mathfrak{p}_\theta, \mathfrak{t}_\theta, \mathfrak{g}_\theta, \mathfrak{r}_\theta$ at uncountable regular cardinals $\theta$.
Motivated by a theorem of Raghavan-Shelah who proved that $\mathfrak{s}_\theta\leq\mathfrak{b}_\theta$, we explore in the first part of the paper the consistency of inequalities comparing $\mathfrak{s}_\theta$ with $\mathfrak{p}_\theta$ and $\mathfrak{g}_\theta$.
In the second part of the paper we study variations of the extender-based Radin forcing to establish several consistency results concerning $\mathfrak{r}_\theta,\mathfrak{s}_\theta$ from hyper-measurability assumptions, results which were previously known to be consistent only from supercompactness assumptions.
In doing so, we answer questions from \cite{MR3436372}, \cite{RagSh1} and \cite{MR3509813}, and improve the large cardinal strength assumptions for results from \cite{MR1632081} and \cite{MR3564375}.
\end{abstract}

\maketitle

\newpage

\section{Introduction}

The study of cardinal characteristics of the continuum has been a prominent subject in modern set theory for many years.
Starting from the 1990s, a new line of research concerning generalized cardinal characteristics associated with uncountable cardinals $\theta$ has started gaining momentum.
The seminal articles \cite{MR1355135} and \cite{MR1251349} are among the first studies in this direction.
This line of research has flourished in recent years, with many advancements made to both the ZFC theory of generalized cardinal characteristics, as well as to consistency results and the related forcing and large cardinal theories.

The purpose of this paper is to study consistency results concerning several cardinal characteristics at a regular uncountable cardinal $\theta$, using methods of forcing with large cardinals.
Especially, we are interested in the invariants $\mathfrak{s}_\theta, \mathfrak{p}_\theta, \mathfrak{g}_\theta, \mathfrak{r}_\theta, \mathfrak{u}_\theta$, and questions of the form: (i) which inequalities between different invariants are provable in ZFC; and (ii) which inequalities can be shown to be consistent, and from which large cardinal assumptions.
Our results are mainly motivated by and build on the recent studies: \cite{RagSh1}, \cite{MR3436372}, \cite{MR3564375} and \cite{MR3509813}.

In the first part of the paper we prove that if $\theta$ is supercompact then one can force $\mathfrak{p}_\theta,\mathfrak{g}_\theta<\mathfrak{s}_\theta<2^\theta$.
In particular, these results show that the bound $\mathfrak{s}_\theta\leq\mathfrak{b}_\theta$, proved by Raghavan-Shelah in \cite{MR3615051}, does not extend to $\mathfrak{s}_\theta$ and the cardinal characteristics $\mathfrak{p}_\theta,\mathfrak{g}_\theta$.
In the second part of the paper we improve the known upper-bounds on the consistency strength of the inequality $\mathfrak{r}_\theta<2^\theta=\theta^{++}$ where $\theta$ is weakly compact but not measurable or measurable but not supercompact.
The former is shown to be consistent from an assumption slightly weaker than the existence of a measurable cardinal $\theta$ with $o(\theta)=\theta^{+3}$, and the latter from $o(\theta)=\lambda$ where $\lambda>\theta$ is a weakly compact cardinal.
The results here build on the extender-based Radin forcing of Merimovich, \cite{MR2783980}, and on ideas from a previous work by the first author and Gitik, \cite{MR3436372}, on the generalized splitting number $\mathfrak{s}_\theta$.

The rest of the paper is organized as follows.
In Section 1 we establish the consistency of the inequalities concerning the generalized cardinal characteristics $\mathfrak{s}_\theta$ and $\mathfrak{p}_\theta,\mathfrak{g}_\theta$.
Section 2 surveys the basics of the extender-based Radin forcing using the ideas and terminology of Merimovich, \cite{MR2783980}.
In Section 3 we use the extender-based Radin forcing to deal with the consistency strength of a small value of the reaping number at weakly compact but not measurable cardinals.
In the last section we deal with the consistency strength of $\mathfrak{r}_\theta<2^\theta=\theta^{++}$ and $\mathfrak{s}_\theta>\theta^+$ at measurable but not supercompact cardinals.

Our notation is mostly standard.
We shall use the Jerusalem forcing notation, so if $p\leq_{\mathbb{P}}q$ then $p$ is weaker than $q$.
The symbol $\lambda\cdot\kappa$ refers to \emph{ordinal multiplication}.
If $A,B\subseteq\theta$ then $A\subseteq^* B$ iff $|A-B|<\theta$.
For a general background in cardinal characteristics we suggest \cite{MR2768685}.
We direct the reader to \cite{MR2768695} regarding Prikry-type forcings in general, and to \cite{MR2783980} for more specific account of extender-based Radin forcing.
We thank the referee for many helpful suggestions which were integrated within the manuscript. We would also like to thank Moti Gitik for valuable discussions on the topic of generalized characteristics, and in particular for helping us to detect an error in one of the theorems which appeared in the original version of the paper.

\newpage

\section{Navigating with gps}

By Raghavan-Shelah \cite{MR3615051}, if $\theta$ is regular and uncountable then $\mathfrak{s}_\theta\leq\mathfrak{b}_\theta$.
Motivated by this result, we explore in this section the possible consistency of similar inequalities between $\mathfrak{s}_\theta$ and the cardinal characteristics $\mathfrak{p}_\theta,\mathfrak{t}_\theta, \mathfrak{g}_\theta$.
We commence with a few definitions:

\begin{definition}
\label{defpkappa} The pseudointersection number. \newline
Let $\theta$ be a regular cardinal.
\begin{enumerate}
\item [$(\aleph)$] A family $\mathscr{F}\subseteq[\theta]^\theta$ has the strong intersection property iff $|\bigcap a|=\theta$ whenever $a\in[\mathscr{F}]^{<\theta}$.
\item [$(\beth)$] A set $S\in[\theta]^\theta$ is a pseudointersection of $\mathscr{F}\subseteq[\theta]^\theta$ iff $S\subseteq^* A$ for every $A\in\mathscr{F}$.
\item [$(\gimel)$] A $\mathfrak{p}_\theta$-family is a collection of sets $\mathscr{F}\subseteq[\theta]^\theta$ with the strong intersection property yet no pseudointersection.
\item [$(\daleth)$] $\mathfrak{p}_\theta$ is the minimal size of a $\mathfrak{p}_\theta$-family.
\end{enumerate}
\end{definition}

A similar definition gives the following concept:

\begin{definition}
\label{deftkappa} The towering number. \newline
Let $\theta$ be a regular cardinal.
\begin{enumerate}
\item [$(\aleph)$] A family $\mathscr{T} = \{T_\alpha:\alpha<\lambda\} \subseteq[\theta]^\theta$ is a $\theta$-tower iff $\alpha<\beta<\lambda \Rightarrow T_\beta\subseteq^* T_\alpha, \mathscr{T}$ has the strong intersection property and no pseudointersection.
\item [$(\beth)$] $\mathfrak{t}_\theta$ is the minimal size of a $\theta$-tower.
\end{enumerate}
\end{definition}

Clearly, $\mathfrak{p}_\theta\leq\mathfrak{t}_\theta$.
It is unknown whether a strict inequality is consistent for some $\theta$.
Both characteristics are \emph{small} in the sense that they are bounded by most cardinal characteristics.
In order to define a dividing line between small and large cardinal characteristics, we suggest a criterion which depends on the cofinality of $2^\theta$.
We shall say that a cardinal characteristic over $\theta$ is small iff it is bounded by $\cf(2^\theta)$.
Actually, we can define the following:

\begin{definition}
\label{defkonig} The K\"onig number.
\begin{enumerate}
\item [$(\aleph)$] The K\"onig number $\mathfrak{k}$ is the first cardinal $\lambda$ such that $2^\lambda>2^\omega$.
\item [$(\beth)$] If $\theta$ is an infinite cardinal then $\mathfrak{k}_\theta$ is the first cardinal $\lambda$ such that $2^\lambda>2^\theta$.
\end{enumerate}
\end{definition}

It is known that $\mathfrak{t}_\theta\leq\mathfrak{k}_\theta$ whenever $\theta=\cf(\theta)$.
The proof is easier when $\theta=\aleph_0$ (see \cite[Theorem 6.14]{MR2768685}), but it holds at uncountable cardinals as well (see \cite{MR1997605}).
Our central definition in this section is the following cardinal characteristic:

\begin{definition}
\label{defskappa} The splitting number. \newline
Let $\theta$ be an infinite cardinal.
\begin{enumerate}
\item [$(\aleph)$] For $B\in[\theta]^\theta$ and $S\subseteq\theta$, we shall say that $S$ splits $B$ if $|S\cap B|=|(\theta-S)\cap B|=\theta$.
\item [$(\beth)$] $\{S_\alpha:\alpha<\kappa\}$ is a splitting family in $\theta$ iff for every $B\in[\theta]^\theta$ there exists an ordinal $\alpha<\kappa$ so that $S_\alpha$ splits $B$.
\item [$(\gimel)$] The splitting number $\mathfrak{s}_\theta$ is the minimal cardinality of a splitting family in $\theta$.
\end{enumerate}
\end{definition}

If $\theta$ is supercompact then one can force $\mathfrak{s}_\theta=\lambda$ for every prescribed $\lambda=\cf(\lambda)>\theta$.
We shall see that this can be done along with $\mathfrak{t}_\theta=\kappa$ (or $\mathfrak{p}_\theta=\kappa$) when $\kappa=\cf(\kappa)<\lambda$.
Let us describe shortly the idea.
The natural way to increase $\mathfrak{s}_\theta$ is the generalized Mathias forcing, as done by Kamo in \cite{Kamo}.
If we apply this method without special care then $\mathfrak{p}_\theta=\mathfrak{s}_\theta$ in the generic extension, since basically each component of Mathias forcing destroys more and more small $\mathfrak{p}_\theta$-families.

In order to maintain some small $\mathfrak{p}_\theta$-families or towers in the generic extension we need a careful bookkeeping.
A different approach is to use the fact that $\mathfrak{t}_\theta\leq\mathfrak{k}_\theta$.
The idea is to choose $\kappa=\cf(\kappa)\geq\theta^+$ and to force $2^\kappa>2^\theta$.
In this way we obtain the upper bound $\kappa\geq\mathfrak{t}_\theta$.
Now we force $\mathfrak{s}_\theta=\lambda>\kappa$ and the result follows.
The same method has been exploited in \cite{MR1997605} in order to force $\mathfrak{t}_\theta<\mathfrak{b}_\theta$.
Unlike the splitting number, there is no need of large cardinals for the inequality $\mathfrak{t}_\theta<\mathfrak{b}_\theta$.

Within the proof of the main theorem we shall use instances of forcing notions which destroy $\mathfrak{p}_\theta$-families (or $\theta$-towers) and forcing notions which add generalized Mathias reals and thus destroy splitting families.
Let us define these components.

\begin{definition}
\label{defbabel} Babel forcing. \newline
Let $\theta$ be a regular cardinal, $\eta=\cf(\eta)>\theta$. \newline
Let $x = \{x_\gamma:\gamma<\eta\}$ be a $\mathfrak{p}_\theta$-family. \newline
A condition $p$ in Babel forcing $\mathbb{B}_x$ is a pair $(a,s)=(a^p,s^p)$ such that $a\in[\theta]^{<\theta}$ and $s\in[\eta]^{<\theta}$.
If $p,q\in\mathbb{B}_x$ then $p\leq q$ iff $a^p\subseteq a^q, s^p\subseteq s^q$ and if $\delta\in a^q-a^p$ then $\delta\in\bigcap\{x_\gamma:\gamma\in s^p\}$.
\end{definition}

It is routine to verify that $\mathbb{B}_x$ is $\theta$-directed-closed.
If $\theta = \theta^{<\theta}$ then $\mathbb{B}_x$ is $\theta^+$-cc.
If $G\subseteq\mathbb{B}_x$ is generic then $x_G = \bigcup\{a:\exists s, (a,s)\in G\}$.
By density arguments, $|x_G|=\theta$ and $x_G\subseteq^* x_\gamma$ for every $\gamma<\eta$.
So if $x$ is a $\mathfrak{p}_\theta$-family (or a tower) then $\mathbb{B}_x$ destroys this property.
If $\theta$ is an indestructible supercompact cardinal then it remains supercompact after forcing with $\mathbb{B}_x$.

\begin{definition}
\label{defmathias} Generalized Mathias forcing. \newline
Let $\theta$ be a measurable cardinal, $\mathcal{U}$ a $\theta$-complete ultrafilter over $\theta$. \newline
A condition $p\in\mathbb{M}_{\mathcal{U}}$ is a pair $(a,A)=(a^p,A^p)$ such that $a\in[\theta]^{<\theta}$ and $A\in\mathcal{U}$.
If $p,q\in\mathbb{M}_{\mathcal{U}}$ then $p\leq q$ iff $a^p\subseteq a^q, A^p\supseteq A^q$ and $a^q-a^p\subseteq A^p$.
\end{definition}

As in the former definition, $\mathbb{M}_{\mathcal{U}}$ is $\theta^+$-cc and $\theta$-directed-closed.
If $G\subseteq\mathbb{M}_{\mathcal{U}}$ is generic then we let $y_G = \bigcup\{a:\exists A, (a,A)\in G\}$.
It can be shown that $y_G\subseteq^* y\vee y_G\subseteq^*(\theta-y)$ for every $y\in[\theta]^\theta\cap V$.
By iterating $\mathbb{M}_{\mathcal{U}}$ over a supercompact cardinal $\theta$ (with $<\theta$-support) one increases $\mathfrak{s}_\theta$.
Of course, the forcing is rendered over a Laver-indestructible supercompact cardinal so that supercompactness is preserved at any stage of the iteration.

If $\lambda=\cf(\lambda)>\theta$ and the iteration is of length $\lambda$ then $\mathfrak{s}_\theta\geq\lambda$ in the generic extension.
Indeed, a family $\mathscr{F}\subseteq[\theta]^\theta$ of size less than $\lambda$ will appear at a bounded stage of the iteration, and the generic $y_G$ added at this stage is not split by $\mathscr{F}$.
Likewise, $\mathfrak{s}_\theta\leq\lambda$ since $\theta$-Cohen subsets are introduced at $\lambda$-many stages of the iteration, and they form a splitting family.
Hence iterating $\mathbb{M}_{\mathcal{U}}$ is a convenient way to set the size of $\mathfrak{s}_\theta$.

The main result of this section reads as follows:

\begin{theorem}
\label{thmpts} Let $\theta$ be a supercompact cardinal. \newline
Assume that $\theta^+\leq\kappa=\cf(\kappa)\leq\lambda=\cf(\lambda)\leq\mu$ and $\cf(\mu)>\theta$. \newline
Then one can force $\mathfrak{p}_\theta=\mathfrak{t}_\theta=\kappa, \mathfrak{s}_\theta=\lambda$ and $2^\theta=\mu$.
\end{theorem}

\par\noindent\emph{Proof}. \newline
We commence with Laver's preparatory forcing which makes $\theta$ indestructible under $\theta$-directed-closed forcing notions, and we may add GCH above $\theta$.
Let $V$ be the universe after this preparation.

Let $\mathbb{C}$ be Cohen forcing for adding $\mu^+$ subsets of $\kappa$.
A condition $p\in\mathbb{C}$ is a partial function from $\mu^+\times\kappa$ into $\{0,1\}$ such that $|{\rm dom}(p)|<\kappa$.
If $p,q\in\mathbb{C}$ then $p\leq_{\mathbb{C}}q$ iff $p\subseteq q$.
Observe that $\mathbb{C}$ is $\kappa$-closed since $\kappa$ is regular.
Likewise, $\mathbb{C}$ is $\kappa^+$-cc by a $\Delta$-system argument and the fact that $\kappa=\kappa^{<\kappa}$.
Consequently, if $G\subseteq\mathbb{C}$ is $V$-generic then $V[G]\models 2^\kappa=\mu^+$.
Moreover, $2^\partial=\partial^+$ for every $\partial\in[\theta,\kappa)$ in $V[G]$, since $\mathbb{C}$ is $\kappa$-closed and hence adds no bounded subsets of $\kappa$.

In $V[G]$ we define a $<\theta$-support iteration $\mathbb{P} = \langle\mathbb{P}_\alpha,\name{\mathbb{Q}}_\beta: \alpha\leq\mu\cdot\lambda, \beta<\mu\cdot\lambda\rangle$ as follows.
Let $\Gamma = \mu\cdot\lambda - \{\mu\cdot\varepsilon:\varepsilon<\lambda\}$.
Let $f:\Gamma\rightarrow(\mu\cdot\lambda)\times\mu$ be onto, such that the following proviso is satisfied:
$$
\forall\beta, f(\beta)=(\gamma,\delta) \Rightarrow \gamma\leq\beta.
$$
This requirement makes sure that $V[G][H_\gamma]\subseteq V[G][H_\beta]$ when we come to choose the $V[G]$-generic set $H$ for the iteration.

By induction on $\beta<\mu\cdot\lambda$ we choose $\mathbb{P}_\beta$ names for posets by the following procedure.
If $\beta=0$ then $\name{\mathbb{Q}}_\beta$ is the name of the empty forcing.
If $\beta\in\Gamma$ and $\beta>0$ then let $\{x^\beta_\zeta:\zeta<\mu\}$ be an enumeration (possibly with repetitions) of all names of $\mathfrak{p}_\theta$-families of size $\eta$ for some $\eta=\cf(\eta)$ so that $\theta^+\leq\eta<\kappa$.
Apply $f(\beta)$ to get a pair $(\gamma,\delta)$, so $\gamma\leq\beta$.
Let $x$ be $x^\gamma_\delta$ and let $\name{\mathbb{Q}}_\beta$ be (a name of) the forcing $\mathbb{B}_x$.
If $\beta\in\mu\cdot\lambda - \Gamma$ then let $\name{\mathcal{U}}_\beta$ be a name of a $\theta$-complete ultrafilter over $\theta$ and let $\name{\mathbb{Q}}_\beta$ be $\mathbb{M}_{\name{\mathcal{U}}_\beta}$.
By the properties of the components of $\mathbb{P}$ one can see that $\mathbb{P}$ is $\theta^+$-complete and $\theta$-directed-closed, so all cardinals are preserved upon forcing with $\mathbb{P}$ and $\theta$ remains supercompact.

Fix a $V[G]$-generic set $H\subseteq\mathbb{P}$.
By a standard argument, our choice of support and posets guarantee that the iteration $\mathbb{P}$ satisfies $\theta^+$-cc.
Since the cofinality of its length is $\lambda$, it follows that every family $\mathcal{F}\subseteq[\theta]^\theta$ of size less than $\lambda$ appears in some intermediate generic extension of the form $V[G][H_\beta]$ of $V[G][H]$, where $H_\beta=H\cap\mathbb{P}_\beta$ and $\beta<\mu\cdot\lambda$.

By the properties of Babel and Mathias forcing we see that $\mathfrak{p}_\theta\geq\kappa$, and similarly $\mathfrak{t}_\theta\geq\kappa$.
Indeed, any $\mathfrak{p}_\theta$-family (including the towers) of size $\eta\in[\theta^+,\kappa)$ has been destroyed by some $\mathbb{B}_x$ along the iteration.
On the other hand, $\mathfrak{p}_\theta,\mathfrak{t}_\theta\leq\kappa$ since $2^\kappa>2^\theta$ in $V[G][H]$ and due to \cite{MR1997605}.
The generalized Mathias components iterated cofinally in the ordinal $\mu\cdot\lambda$ set $\mathfrak{s}_\theta=\lambda$, so we are done.

\hfill \qedref{thmpts}

We conclude this section with the consistency of $\mathfrak{g}_\theta<\mathfrak{s}_\theta$.
This will be proved, again, under the assumption that $\theta$ is supercompact.
Let us begin with the definition of $\mathfrak{g}_\theta$, the generalized groupwise density:

\begin{definition}
\label{defg} Let $\theta$ be a regular uncountable cardinal.
\begin{enumerate}
\item [$(\aleph)$] A family $\mathcal{G}\subseteq[\theta]^\theta$ is groupwise dense iff it is downward closed under $\subseteq^*$ and for every increasing $f\in{}^\theta\theta$ there exists $B\in[\theta]^\theta$ such that $\bigcup\{[f(\alpha),f(\alpha+1)):\alpha\in B\}\in\mathcal{G}$.
\item [$(\beth)$] The groupwise density number $\mathfrak{g}_\theta$ is the smallest cardinal $\kappa$ such that there is a collection $(\mathcal{G}_\gamma:\gamma\in\kappa)$ of groupwise dense families in $[\theta]^\theta$ with empty intersection.
\end{enumerate}
\end{definition}

Lest $\theta=\aleph_0$ one can force $\mathfrak{b}<\mathfrak{g}$ and even $\mathfrak{u}<\mathfrak{g}$, so in some sense $\mathfrak{g}$ is quite large.
On the other hand, $\mathfrak{g}\leq\cf(2^\omega)$ so in some sense it is small.
The generalization of the latter fact yields the consistency of $\mathfrak{g}_\theta<\mathfrak{s}_\theta$.

\begin{theorem}
\label{thmgs} Let $\theta$ be a supercompact cardinal. \newline
Then one can force $\mathfrak{g}_\theta<\mathfrak{s}_\theta$, and the gap can be arbitrarily large.
\end{theorem}

\par\noindent\emph{Proof}. \newline
Firstly we show that $\mathfrak{g}_\theta\leq\cf(2^\theta)$.
This statement is proved exactly as the parallel statement for $\theta=\aleph_0$, see \cite{MR2768685}.
We indicate that it holds at every regular cardinal $\theta$, and supercompactness is not needed at this stage.
To prove this claim, suppose that $\mathcal{F}\subseteq[\theta]^\theta$ and $|\mathcal{F}|<2^\theta$.
We shall construct a groupwise dense family $\mathcal{G}$ such that $\mathcal{G}\cap\mathcal{F}=\varnothing$.
We set:
$$
\mathcal{G} = \{x\in[\theta]^\theta:\forall y\in\mathcal{F}, \neg(y\subseteq^* x)\}.
$$
\begin{enumerate}
\item [$(a)$] $\mathcal{G}$ is $\subseteq^*$-downward closed since $\subseteq^*$ is transitive.
\item [$(b)$] $\mathcal{G}\cap\mathcal{F}=\varnothing$ since $\subseteq^*$ is reflexive.
\item [$(c)$] $\mathcal{G}$ is groupwise dense.
For this, let $f\in{}^\theta\theta$ be an increasing function, and let $\mathcal{A}\subseteq[\theta]^\theta$ be an almost disjoint family of size $2^\theta$.
For each $A\in\mathcal{A}$ let $A_f = \bigcup\{[f(\alpha),f(\alpha+1)):\alpha\in A\}$. Notice that $A,B\in\mathcal{A}\Rightarrow |A_f\cap B_f|<\theta$.
We claim that for some $A\in\mathcal{A}$ the statement $\forall y\in\mathcal{F}, \neg(y\subseteq^* A_f)$ holds true.
If not, there are distinct $A,B\in\mathcal{A}$ such that for some $y\in\mathcal{F}$ we have $y\subseteq^* A_f \wedge y\subseteq^* B_f$, since $|\mathcal{F}|<2^\theta=|\mathcal{A}|$.
But then $|A_f\cap B_f|=\theta$, which is impossible.
So choose $A\in\mathcal{A}$ such that $\forall y\in\mathcal{F}, \neg(y\subseteq^* A_f)$.
This means that $A_f\in\mathcal{G}$, thus showing that $\mathcal{G}$ is groupwise dense.
\end{enumerate}
Express $[\theta]^\theta$ as $\bigcup\{\mathcal{F}_\beta:\beta<\cf(2^\theta)\}$ such that $|\mathcal{F}_\beta|<2^\theta$ for each $\beta<\cf(2^\theta)$.
By the above consideration, for every $\beta<\cf(2^\theta)$ choose a groupwise dense family $\mathcal{G}_\beta$ such that $\mathcal{F}_\beta\cap\mathcal{G}_\beta=\varnothing$.
Notice that $\bigcap\{\mathcal{G}_\beta:\beta<\cf(2^\theta)\}=\varnothing$ since the union of all the $\mathcal{F}_\beta$ covers $[\theta]^\theta$.
It follows that $\mathfrak{g}_\theta\leq\cf(2^\theta)$, as desired.

Secondly, choose $\kappa=\cf(\kappa)>\theta$ and $\mu>\cf(\mu)=\kappa$.
At this stage we need the assumption that $\theta$ is supercompact.
By Theorem \ref{thmpts} we can force $2^\theta=\mu$ and $\mathfrak{s}_\theta>\kappa$.
Since $\mathfrak{g}_\theta\leq\kappa$, the inequality $\mathfrak{g}_\theta<\mathfrak{s}_\theta$ obtains.

\hfill \qedref{thmgs}

\begin{remark}
\label{r0} We note that $\mathfrak{g}_\theta<\mathfrak{s}_\theta$ and $\mathfrak{p}_\theta<\mathfrak{s}_\theta$ are consistent with $\cf(2^\theta)=2^\theta$.
For example, collapsing $2^\theta$ to $\mathfrak{s}_\theta^+$ in the above theorem preserves these inequalities, a fact which follows from the completeness of the collapse.
\end{remark}

The results of this section are proved under the assumption that there exists a supercompact cardinal in the ground model.
It seems, however, that much less is required.
Basically, we wish to increase $\mathfrak{s}_\theta$ to $\theta^{++}$ while forcing over a universe in which $2^\theta=\mu$ and $\mu>\cf(\mu)=\theta^+$.
The methods of the next sections suggest that a measurable cardinal with sufficiently large Mitchell order will suffice.

The consistency strength of the inequalities forced in the present section is at least a measurable cardinal $\kappa$ with $o(\kappa)=\kappa^{++}$, since $\mathfrak{p}_\theta,\mathfrak{g}_\theta<\mathfrak{s}_\theta$ implies that $\mathfrak{s}_\theta\geq\theta^{++}$.
The consistency strength of $\mathfrak{s}_\theta=\theta^{++}$ is exactly $o(\kappa)=\kappa^{++}$, but here we force when $2^\theta$ is a singular cardinal $\mu$, and increasing $\mathfrak{s}_\theta$ in this environment requires a bit more.

\begin{question}
\label{qgs} What is the consistency strength of the inequality $\mathfrak{g}_\theta<\mathfrak{s}_\theta$?
\end{question}

\newpage

\section{Extender-based Radin forcing}

In this section we give a short account of the extender-based Radin forcing, using the argot of Carmi Merimovich.
Apart of setting our terminology we shall explain how to get the basic properties of this forcing notion under the assumption of $2^\kappa=\kappa^{++}$ in the ground model.
This will be needed for our main theorem, and it differs from \cite{MR2783980} in which GCH is assumed in the ground model.

Radin forcing has been published in \cite{MR670992}, and became a central tool for proving combinatorial statements at small large cardinals.
Extender-based Radin forcing implements the ideas of the extender-based Prikry forcing on Radin forcing.
We shall try to describe the basic idea and en route fix notation and prove some simple facts.
All the results in this section are either due to Merimovich, \cite{MR2783980}, or a slight modification of his results. \newline

\noindent \textbf{Extenders.}

An extender is a directed system of ultrafilters and embeddings.
Our extenders come from elementary embeddings.
Suppose that $\kappa$ is measurable, $\jmath:V\rightarrow M$ is elementary and $\kappa={\rm crit}(\jmath)$.
We shall describe the extender derived from this embedding.
For every $a\in M$ one can define a $\kappa$-complete ultrafilter $E(a)$ over $\kappa$ by letting $A\in E(a)$ iff $a\in\jmath(A)$.
In particular, if $\alpha$ is an ordinal below $\jmath(\kappa)$ then $E(\alpha)$ is such an ultrafilter.
If $\alpha\in\kappa, E(\alpha)$ would be principal, so we focus on ordinals $\alpha\in[\kappa,\jmath(\kappa))$.
In the case of $\alpha=\kappa, E(\alpha)$ is normal.

The extender $E$ derived from $\jmath$ consists of the ultrafilters $E(\alpha)$ for $\alpha\in[\kappa,\jmath(\kappa))$.
It also contains a collection of natural embeddings between the ultrafilters.
Suppose that $\alpha,\beta\in[\kappa,\jmath(\kappa))$.
We shall say that $E(\alpha)\leq^{\rm RK}E(\beta)$ iff there is a projection $\pi=\pi_{\beta\alpha}:\kappa\rightarrow\kappa$ such that $\jmath(\pi)(\beta)=\alpha$.
The symbol $\leq^{\rm RK}$ stands for Rudin-Keisler.
If $\alpha\in\jmath(\kappa)$ then $\alpha$ is called a generator of $E$ iff there is no $\beta<\alpha$ so that $E(\alpha)\leq^{\rm RK}E(\beta)$.
A good example is $\kappa$, the critical point of $\jmath$.

An extender $E$ is the set of ultrafilters $E(\alpha)$ for $\alpha\in[\kappa,\jmath(\kappa))$ together with the projections $\pi_{\beta\alpha}$ for every $\alpha,\beta\in[\kappa,\jmath(\kappa))$ such that $E(\alpha)\leq^{\rm RK}E(\beta)$.
This collection of ultrafilters and embeddings forms a directed system.
Let $M_E$ be the direct limit of the system.
For each ordinal $\alpha\in[\kappa,\jmath(\kappa))$ let $\jmath_{E(\alpha)}:V\rightarrow M_{E(\alpha)}\cong{\rm Ult}(V,E(\alpha))$ be the canonical embedding derived from $E(\alpha)$, and let $\jmath_E:V\rightarrow M_E$ be the corresponding embedding of the direct limit.
A natural elementary map $k:M_E\rightarrow M$ is defined by $k(\jmath_Ef(\alpha))=\jmath f(\alpha)$.
Notice that $k$ is applied here to the equivalence class of some function $f\in{}^\kappa\kappa$.

The critical point of an extender $E$ is the first ordinal $\delta$ such that $\delta<\jmath_E(\delta)$.
It will be denoted by ${\rm crit}(E)$ or ${\rm crit}(\jmath_E)$.
The height of an extender $E$ is defined as $\sup\{\alpha\in\jmath(\kappa):\alpha$ is a generator of $E\}$.
It will be denoted by $\sigma(E)$.
We call $E$ a $(\kappa,\lambda)$-extender iff $\kappa={\rm crit}(\jmath_E)$ and $\lambda=\sigma(E)$.
We shall say that $E$ is \emph{short} iff $\sigma(E)\leq\jmath_E(\kappa)$.

An important preorder defined on ultrafilters is the so-called Mitchell order.
It can be applied to extenders as well.
suppose that $E$ and $F$ are two extenders.
We shall say that $E\lhd F$ iff $E\in M_F$, where $M_F$ is the transitive collapse of ${\rm Ult}(V,F)$.
We shall use $\lhd$-increasing sequences of extenders over a measurable cardinal $\kappa$.
The length of such a sequence calibrates the strength of our assumption on $\kappa$.
We say that $o(\kappa)\geq\lambda$ iff there exists a $\lhd$-increasing sequence $(E_\eta:\eta\in\lambda)$ of extenders such that the critical point of the extenders is $\kappa$ and $(\sigma(E_\eta):\eta\in\lambda)$ is strictly increasing and $\sigma(E_\eta)<\lambda$ for every $\eta\in\lambda$.

Our aim is to begin with a measurable cardinal $\kappa$ such that $o(\kappa)=\kappa^{+3}$ in the ground model $V$.
From this assumption we shall force $\mathfrak{r}_\kappa=\kappa^+<\kappa^{++}=2^\kappa$ while preserving the fact that $o(\kappa)=\kappa^{+3}$.
Let $(e_\eta:\eta\in\kappa^{+3})$ be a sequence of extenders over $\kappa$ which exemplifies $o(\kappa)=\kappa^{+3}$.
By cutting off an initial segment and shrinking further if needed we may assume that $\sigma(e_\eta)>\kappa^{++}$ and $\eta\leq\sigma(e_\eta)$ for every $\eta\in\kappa^{+3}$.

Classical constructions enable us to force over $V$ and obtain $2^\kappa=\kappa^{++}$ in $V[G]$ while making sure that the embeddings associated with the sequence of extenders can be lifted.
Namely, if $\imath:V\rightarrow M_{e_\eta}$ is the embedding associated with $e_\eta$ then $\imath$ extends to an elementary embedding $\imath^+:V[G]\rightarrow M_{e_\eta}[H]$ where $H=\imath^+(G)$.
Such a construction is described in \cite{MR1007865}.
If $a\in M_{e_\eta}$ then one can define the $\imath^+_{e_\eta(a)}$-derived ultrafilter $E_\eta(a)$ by $A\in E_\eta(a)$ iff $a\in\imath^+_{e_\eta(a)}(A)$.
One can verify that $E_\eta(a)$ extends $e_\eta(a)$ from $V$.
Thence it is possible to define the $\imath^+_{e_\eta}$-derived extender which we call $E_\eta$.
It is well-known that if $(e_\eta:\eta\in\kappa^{+3})$ is $\lhd$-increasing then the derived sequence $(E_\eta:\eta\in\kappa^{+3})$ is $\lhd$-increasing as well. \newline

\noindent \textbf{Extender-based Radin forcing.}

The main forcing tool for obtaining our consistency results is the extender-based Radin forcing, introduced by Carmi Merimovich in \cite{MR2783980}.
We briefly describe the forcing and its main properties, and refer the reader to \cite[Section 4]{MR2783980} for a full description and detailed proofs.

For concreteness, we describe the forcing over a measurable cardinal $\kappa$ with $o(\kappa)=\kappa^{+3}$.
The description can be easily generalized to $o(\kappa)=\lambda$ for an arbitrary cardinal $\lambda$.
Assume therefore that $\kappa$ is measurable and $o(\kappa)=\kappa^{+3}$.
This means that there exists a $\lhd$-increasing sequence of extenders over $\kappa$.
Let $\bar{E}=(E_\eta:\eta\in\kappa^{+3})$ be such a sequence over $\kappa$.
Define $\epsilon=\sup\{\jmath_{E_\eta}(\kappa):\eta\in\kappa^{+3}\}$.
For every $\alpha\in[\kappa,\epsilon)$ let $\bar{\alpha}=\langle\alpha\rangle^\frown\langle E_\eta:\eta\in\kappa^{+3}, \alpha<\jmath_{E_\eta}(\kappa)\rangle$.
To each condition in the extender based Radin forcing by $\bar{E}$, there is an associated support set $\bar{d}=\{\bar{\alpha}:\alpha\in d\}$, where $d\in[\kappa^{+3}]^{\leq\kappa}$.
We only consider sets $d$ with $\kappa\in d$.

We introduce some additional notation to describe the conditions of the forcing.
For every $\bar{\alpha}$ and for each $\eta\in[0,\kappa^{+3})$ define the initial segement $R_\eta(\bar{\alpha})$ of $\bar{\alpha}$ as $\langle\alpha\rangle^\frown\langle E_\zeta:\zeta\in\eta,\alpha<\sigma(E_\zeta)\rangle$.
Then, for $\eta\in\kappa^{+3}$ and $d\in[\kappa^{+3}]^{\leq\kappa}$ we define the maximal coordinate of $d$ (with respect to $\eta$) to be the following function:
$$
{\rm mc}_\eta(d)=\{\langle\jmath_{E_\eta}(\bar{\alpha}), R_\eta(\bar{\alpha})\rangle: \bar{\alpha}\in d,\alpha<\jmath_{E_\eta}(\kappa)\}.
$$
Notice that $R_\eta(\bar{\alpha})\in M_{E_\eta}$ and moreover ${\rm mc}_\eta(d)\in M_{E_\eta}$.

In essence, the function ${\rm mc}_\eta(d)$ is an augmentation of the function $(\jmath_{E_\eta}\upharpoonright d)^{-1}$ which maps an ordinal of the form $\jmath_{E_\eta}(\alpha)$ back to $\alpha$, for $\alpha\in d$.
The main difference is that $(\jmath_{E_\eta}\upharpoonright d)^{-1}$ maps ordinals to ordinals, while the function ${\rm mc}_\eta(d)$ further indicates an extender sequence $R_\eta(\bar{\alpha})$ which is paired with $\alpha$ in the image of the function.

Next, \cite[Definition 4.3]{MR2783980} defines a set ${\rm OB}(d)\subseteq V_\kappa$ of partial functions $\nu:\bar{d}\rightarrow V_\kappa$, which reflect the basic properties of the functions ${\rm mc}_\eta(d)\in M_{E_\eta}$ for various $\eta\in\kappa^{+3}$.
We only specify here some of these properties, and refer the reader to the full definition in \cite{MR2783980}.
We first set $\mathfrak{R}=\{\bar{\tau}=\langle\tau\rangle^\frown\bar{e}\in V_\kappa:\bar{e}$ is an extender sequence with critical point $\mu$ and $\tau\in\mu^{+3}\}$.
For $\bar{\tau}=\langle\tau\rangle^\frown\bar{e}\in\mathfrak{R}$ we denote the ordinal $\tau$ by $\bar{\tau}_0$.

The set ${\rm OB}(d)$ consists of partial functions $\nu:\bar{d}\rightarrow\mathfrak{R}$ so that for every $\bar{\alpha}\in{\rm dom}(\nu)$ we have $\nu(\bar{\alpha})=\bar{\tau}\in\mathfrak{R}$.
One also requires that $\bar{\kappa}\in{\rm dom}(\nu)$ and $\nu(\bar{\kappa})=\bar{\mu}$ where $\bar{\mu}_0$ is the critical point of the extenders in $\bar{e}$.
Finally, one requires that $|\nu|\leq\bar{\mu}_0$ and for each $\bar{\alpha}\in{\rm dom}(\nu)$, if $\nu(\bar{\alpha})=\langle\tau\rangle^\frown\bar{e}'$ then $\tau<\bar{\mu}_0^{+3}$ and $\bar{e}'$ is an end segment of $\bar{e}$.

It is straightforward to verify that ${\rm mc}_\eta(d)\in \jmath_{E_\eta}({\rm OB}(d))$ for all $\eta\in\kappa^{+3}$.
The object ${\rm mc}_\eta(d)$ is taken to be the seed of the measure $E_\eta(d)$ on ${\rm OB}(d)$.
Namely, $E_\eta(d)$ consists of sets $A\subseteq{\rm OB}(d)$ so that ${\rm mc}_\eta(d)\in\jmath_{E_\eta}(A)$.
Finally, we define $E(d)=\bigcap\{E_\eta(d):\eta\in\kappa^{+3}\}$.

We are ready to define the extender-based Radin forcing $\mathbb{P}_{\bar{E}}$.
A condition $p\in\mathbb{P}_{\bar{E}}$ consists of a lower part $p_\leftarrow$ and a top element $p_\rightarrow$, so $p={p_\leftarrow}^\frown p_\rightarrow$.
The top element $p_\rightarrow$ is a pair $(f,T)$ where $f:\bar{d}\rightarrow{}^{<\omega}\mathfrak{R}$ for some $d\in[\kappa^{+3}]^{\leq\kappa}$ and $T\subseteq{\rm OB}(d)^{<\omega}$ is a tree associated with $E(d)$.
That is, $T\subseteq{\rm OB}(d)^{<\omega}$ is a tree consisting of finite sequences $t=\langle\nu_0,\ldots,\nu_{\ell-1}\rangle$ of functions in ${\rm OB}(d)$, so that the successor set ${\rm succ}_T(t)=\{\nu:t^\frown\langle\nu\rangle\in T\}$ of each element $t\in T$ belongs to the filter $E(d)$.
We note that if $d'\subseteq d$ then the restricted tree $T\upharpoonright d'=\{\langle\nu_0\upharpoonright\bar{d}',\ldots, \nu_{\ell-1}\upharpoonright\bar{d}'\rangle: \langle\nu_0,\ldots,\nu_{\ell-1}\rangle\in T\}$ satisfies the same requirement with respect to $E(d')$.

The lower part $p_\leftarrow$ is a finite sequence of the from $\langle p_0,\ldots,p_{\ell-1}\rangle\in V_\kappa$.
Each $p_i$ is a pair $(f_i,T_i)$ which looks like the single top-element described above, with respect to some extender-sequence $\bar{e}\in\mathfrak{R}$ over some measurable cardinal below $\kappa$.
The function $f$ which appears in the pair $(f,T)$ will be called \emph{the Cohen part} of the condition.

Suppose that $p,q\in\mathbb{P}_{\bar{E}}$.
We wish to define the orders $\leq$ and $\leq^*$ (the forcing order and the pure order).
Let us begin with the pure order $\leq^*$.
We assume, first, that $p_\leftarrow=q_\leftarrow=\varnothing$, so $p=(f,T)$ and $q=(g,U)$.
We shall say that $p\leq^*q$ iff $f\subseteq g$ and $U\upharpoonright{\rm dom}f\subseteq T$.
So we can extend the domain of the function, shrink the tree, but we do not change the Radin sequences in old coordinates.
Assume now that $p={p_\leftarrow}^\frown p_\rightarrow, q={q_\leftarrow}^\frown q_\rightarrow$ and $p_\leftarrow,q_\leftarrow$ are not necessarily empty.
We define $p_\leftarrow\leq^*q_\leftarrow$ iff $p_\leftarrow=\langle p_0,\ldots,p_{k-1}\rangle, q_\leftarrow=\langle q_0,\ldots,q_{\ell-1}\rangle, k\leq\ell$ and $p_i\leq^*q_i$ for every $i<k$.
Finally, $p\leq^*q$ iff $p_\leftarrow\leq^*q_\leftarrow$ and $p_\rightarrow\leq^*q_\rightarrow$.

On top of the pure extension relation $\leq^*$, the definition of the forcing order $\leq$ is based on the concept of \emph{one-point extension}, given in \cite[Definition 4.5]{MR2783980}.
Again, we only give a rough description here and refer the reader to the definition in \cite{MR2783980} for complete details.

For $\langle\nu\rangle\in T^{p_\rightarrow}$ with $\nu(\bar{\kappa})=\langle\mu_0\rangle^\frown\bar{e}$, the one-point extension $p_{\langle\nu\rangle}$ of $p$ adds, essentially, two new pieces of data to the condition $p$.
The first is a new component $(f_{i'},T_{i'})$ added to $p_{\leftarrow}$, assigned to the extender sequence $\bar{e}$ on $\mu_0$, which is given by $\nu(\bar{\kappa})$.
The second is an update of the function $f$ from the top part $p_\rightarrow=(f,T)$ of $p$, which is rendered as follows.
For each $\bar{\alpha}\in{\rm dom}(\nu)$ we add $\nu(\bar{\alpha})=\langle\tau\rangle^\frown\bar{e}'$ to the finite sequence $f(\alpha)$.
Similarly, one defines a one-point extension $p_{\langle\nu\rangle}$ for $\langle\nu\rangle\in T_i$ in one of the lower components $(f_i,T_i)\in p_\leftarrow$.
Finally, if $p,q\in\mathbb{P}_{\bar{E}}$ then $p\leq q$ iff $q$ is obtained from $p$ by a finite sequence of pure extensions and one-point extensions. \newline

\noindent \textbf{Basic properties.}

Let us give a short summary of some of the basic properties of the extender-based Radin forcing.
We do not make the assumption that $2^\kappa=\kappa^+$ in the ground model, and at this point we differ from \cite{MR2783980}.

\begin{lemma}
\label{lemccc} The forcing notion $\mathbb{P}_{\bar{E}}$ is $\kappa^{++}$-cc even if $2^\kappa>\kappa^+$ in the ground model.
\end{lemma}

\par\noindent\emph{Proof}. \newline
Let $\{p_\alpha:\alpha\in\kappa^{++}\}$ be a subset of $\mathbb{P}_{\bar{E}}$.
The number of possible lower parts of conditions in $\mathbb{P}_{\bar{E}}$ is $\kappa$, so we may assume without loss of generality that if $\alpha<\beta<\kappa^{++}$ then $p_{\alpha\leftarrow}=p_{\beta\leftarrow}$.
Denote the top element $p_{\alpha\rightarrow}$ of $p_\alpha$ by $(f_\alpha,T_\alpha)$ for every $\alpha\in\kappa^{++}$.

The domain of $f_\alpha$ is a set $\bar{d}_\alpha$ of size at most $\kappa$.
By the Delta-system lemma there are $\alpha<\beta<\kappa^{++}$ such that $f_\alpha$ and $f_\beta$ are compatible as functions.
By the definition of $\mathbb{P}_{\bar{E}}$ this means that $p_{\alpha\rightarrow}\parallel p_{\beta\rightarrow}$.
From our assumption that $p_{\alpha\leftarrow}=p_{\beta\leftarrow}$ we infer that $p_\alpha\parallel p_\beta$.

\hfill \qedref{lemccc}

The forcing notion $\mathbb{P}_{\bar{E}}$ preserves $\kappa^+$ and has the Prikry property.
This is true, again, even if $2^\kappa>\kappa^+$ in $V$.
A key point here is the following:

\begin{lemma}
\label{lemproper} Let $\chi$ be a sufficiently large regular cardinal above $\kappa$. \newline
For every $N\prec\mathcal{H}(\chi)$ such that $|N|=\kappa,{}^{<\kappa}N\subseteq N$ and $\bar{E},\mathbb{P}_{\bar{E}}\in N$, if $p\in\mathbb{P}_{\bar{E}}\cap N$ then there exists a pure extension $q$ of $p$ such that $q$ is $(N,\mathbb{P}_{\bar{E}})$-generic.
\end{lemma}

\hfill \qedref{lemproper}

This lemma gives a kind of properness which serves both for proving the Prikry property and for showing that $\kappa^+$ is preserved.
For the proof of this lemma we refer the reader to \cite[Section 4]{MR2783980}.
An examination of the proof shows that the value of $2^\kappa$ in $V$ plays no role.

Suppose that $p\in\mathbb{P}_{\bar{E}}$ and express $p$ as ${p_{\leftarrow}}^\frown p_\rightarrow$.
Let $\mathbb{P}_{\bar{E}}/p$ denote the set $\{q\in\mathbb{P}_{\bar{E}}:q\geq p\}$.
Observe that $\mathbb{P}_{\bar{E}}/p\cong \mathbb{P}_{\bar{E}}/p_\leftarrow\times\mathbb{P}_{\bar{E}}/p_\rightarrow$ by the mapping $r\mapsto(r_\leftarrow,r_\rightarrow)$.
A finite sequence of applications of this fact shows that if $p=\langle p_0,\ldots,p_{\ell-1}\rangle$ and $0<i<\ell-1$ then $\mathbb{P}_{\bar{E}}/p$ factors to $\mathbb{P}_{\bar{E}}/(p_0,\ldots,p_i)\times \mathbb{P}_{\bar{E}}/(p_{i+1},\ldots,p_{\ell-1})$.

\begin{corollary}
\label{corcard} Forcing with $\mathbb{P}_{\bar{E}}$ preserves cardinals.
\end{corollary}

\par\noindent\emph{Proof}. \newline
If $\lambda\geq\kappa^{++}$ then $\lambda$ is preserved by virtue of Lemma \ref{lemccc}.
If $\lambda=\kappa^+$ one can use Lemma \ref{lemproper}.
Suppose that $\lambda<\kappa$ and choose a condition $p\in\mathbb{P}_{\bar{E}}$ such that the top element of the lower part $p_\leftarrow$ which we will denote by $p_{\leftarrow\rightarrow}$ is above $\lambda$ (that is, for some measurable cardinal $\tau\in(\lambda,\kappa)$ the extender sequence used for $p_{\leftarrow\rightarrow}$ is defined over $\tau$).

Recall that $\mathbb{P}_{\bar{E}}/p$ factors to the product $\mathbb{P}_{\bar{E}}/p_{\leftarrow\rightarrow}\times \mathbb{P}_{\bar{E}}({p_{\leftarrow\rightarrow}}^\frown p_\rightarrow)$.
Now use the facts that $\mathbb{P}_{\bar{E}}({p_{\leftarrow\rightarrow}}^\frown p_\rightarrow)$ is $\lambda$-complete and the Prikry property to conclude that $\lambda$ is preserved.
Finally, $\kappa$ is a limit ordinal and every $\lambda<\kappa$ is preserved, so $\kappa$ is preserved as well.

\hfill \qedref{corcard}

We conclude this section with a few comments.
The first one is that the length of $\bar{E}$ plays no role in the proof of the above lemmata and corollary.
Namely, $\mathbb{P}_{\bar{E}}$ preserves cardinals regardless of the length of $\bar{E}$.
However, $\ell g(\bar{E})$ is very important when computing the cofinality of $\kappa$ in the generic extension.

Let $\delta=\ell g(\bar{E})$.
Merimovich proves in \cite[Section 5]{MR2783980} that $\cf^{V}(\delta)$ determines $\cf^{V[G]}(\kappa)$ in the following way.
If $\delta$ is a successor ordinal then $\cf^{V[G]}(\kappa)=\omega$ as one can extract an $\omega$-cofinal sequence in $\kappa$ from the last element of $\bar{E}$.
If $\delta$ is a limit ordinal and $\cf(\delta)<\kappa$ then $\cf^{V[G]}(\kappa)=\cf(\delta)$.
In the case of $\cf(\delta)=\kappa$ one can construct back again a short cofinal sequence in $\kappa^{V[G]}$ and $\cf^{V[G]}(\kappa)=\omega$.
However, if $\cf(\delta)>\kappa$ then $\kappa$ remains regular in $V[G]$.
The reason is that if $p$ forces that $\name{f}$ is a function from $\zeta$ into $\kappa$ for some $\zeta\in\kappa$ then there is a direct extension of $p$ which forces that the range of $\name{f}$ is bounded in $\kappa$.
Finally, if $\cf(\delta)$ is sufficiently large so that there are repeat points along the sequence $\bar{E}$ then $\kappa$ will be measurable in $V[G]$.
This issue will be dealt with in the last section.

The second comment is about the generic objects added by the extender-based Radin forcing.
Recall that $\epsilon=\sup\{\jmath_{E_\eta}(\kappa):\eta\in\kappa^{+3}\}$ and for each $\alpha\in[\kappa,\epsilon)$ we let $\bar{\alpha}$ be an ordered pair in which the first element is the ordinal $\alpha$.
Given $\bar{\alpha}$ we shall write $\bar{\alpha}_0$ for the ordinal $\alpha$.
With this notation at hand let $G$ be a generic subset of $\mathbb{P}_{\bar{E}}$.
For every $\alpha\in[\kappa,\epsilon)$ let $G^\alpha = \bigcup\{f^{p_\rightarrow}(\bar{\alpha}):p\in G,\bar{\alpha}\in{\rm dom}(f^{p_\rightarrow})\}$.
Bearing in mind that $f^{p_\rightarrow}$ is a function into ${}^{<\omega}\mathfrak{R}$ we see that $f^{p_\rightarrow}(\bar{\alpha})$ points to a finite set of coordinates.
The Radin sequence associated with $\bar{\alpha}$ is defined by:
$$
C^\alpha = \{\bar{\nu}_0:\bar{\nu}\in G^\alpha\}.
$$
These sequences show immediately that $V[G]\models 2^\kappa=|\epsilon|$ since $C^\alpha\neq C^\beta$ whenever $\alpha\neq\beta$ as follows from density arguments.
In the specific case of $\alpha=\kappa$ one can see that $C^\kappa$ is a club subset of $\kappa$ in $V[G]$, and we call it \emph{the Radin club}.

\newpage

\section{Consistency strength for the reaping number}

The main objective of the current section is to give an upper bound on the consistency strength of the statement $\mathfrak{r}_\kappa<2^\kappa$ where $\kappa=\cf(\kappa)>\aleph_0$ is a limit cardinal and $2^\kappa=\kappa^{++}$.
As mentioned in the introduction, a recent result of Raghavan and Shelah produces a model of $\mathfrak{u}_\kappa<2^\kappa$ (and a fortiori $\mathfrak{r}_\kappa<2^\kappa$) from a measurable cardinal.
We shall use a different forcing construction, which will give the above inequality at weakly compact cardinals which are not measurable, and moreover compatible with $2^\kappa=\kappa^{++}$.
We commence with the central definition of this section:

\begin{definition}
\label{defreaping} The reaping number. \newline
Let $\kappa$ be an infinite cardinal.
\begin{enumerate}
\item [$(\aleph)$] $\{R_\alpha:\alpha\in\lambda\}\subseteq[\kappa]^\kappa$ is a reaping family in $\kappa$ iff there is no single $B\in[\kappa]^\kappa$ which splits all the elements of the family.
\item [$(\beth)$] The reaping number $\mathfrak{r}_\kappa$ is the minimal cardinality of a reaping family in $\kappa$.
\end{enumerate}
\end{definition}

The above definition is based on the splitting property, and indeed $\mathfrak{r}_\kappa$ is the dual of $\mathfrak{s}_\kappa$.
Suppose that $\kappa=\cf(\kappa)>\aleph_0$.
It has been proved in \cite{MR3436372} that the consistency strength of the statement $\mathfrak{s}_\kappa>\kappa^+$ is exactly $o(\kappa)=\kappa^{++}$.
The duality between $\mathfrak{s}_\kappa$ and $\mathfrak{r}_\kappa$ suggests that the statement $\mathfrak{r}_\kappa<2^\kappa$ will behave similarly.

Actually, the connection is more palpable.
Recall that $\binom{\lambda}{\kappa}\rightarrow\binom{\lambda}{\kappa}_2$ means that for every $c:\lambda\times\kappa\rightarrow 2$ there are $A\in[\lambda]^\lambda,B\in[\kappa]^\kappa$ such that $c\upharpoonright(A\times B)$ is constant.
In the model of the previous section, if one forces $\mathfrak{s}_\kappa=\kappa^{++}$ and $2^\kappa=\kappa^{+3}$ then necessarily $\mathfrak{r}_\kappa<2^\kappa$.
Indeed, the process of adding $\kappa^{++}$ generalized Mathias reals forces $\binom{2^\kappa}{\kappa}\rightarrow\binom{2^\kappa}{\kappa}_2$, while $\mathfrak{r}_\kappa=2^\kappa$ implies $\binom{2^\kappa}{\kappa}\nrightarrow\binom{2^\kappa}{\kappa}_2$.
Hence we have a model in which $\mathfrak{r}_\kappa<2^\kappa$, and tentatively the consistency strength of this statement is closed to that of $\mathfrak{s}_\kappa=\kappa^{++}$.
The fact that we blow up $2^\kappa$ to a larger value requires a bit more, and the main theorem of this section says that the statement $\mathfrak{r}_\kappa<2^\kappa$ where $\kappa$ is an uncountable strongly inaccessible is consistent relative to the existence of a measurable cardinal $\kappa$ of Mitchell order $o(\kappa)=\kappa^{+3}$.

We indicate that the statement $\mathfrak{r}<2^\omega$ has no consistency strength, and the consistency strength of $\mathfrak{r}_\lambda<2^\lambda$ where $\lambda$ is a strong limit singular cardinal is exactly $o(\kappa)=\kappa^{++}$ as proved in \cite{1156}.
For regular cardinals which are not inaccessible, the inequality $\mathfrak{r}_\lambda<2^\lambda$ is forced from a measurable cardinal in the ground model in \cite{RagSh1}.

Though it is possible to compute the strength of a model in which $\mathfrak{s}_\kappa=\kappa^{++}<2^\kappa$ and employ the above arguments about the polarized relation, we shall use a more direct approach.
We shall force with an extender sequence over a ground model $V$ with some prescribed properties, resulting in a model $V[G]$.
The forcing will add a collection of $\kappa^+$-many generating sets to ultrafilters which appear in intermediate extensions of $V$.
We then show that this collection of generating sets is a reaping family at $\kappa$, thus $\mathfrak{r}_\kappa=\kappa^+<\kappa^{++}=2^\kappa$ in $V[G]$.

We can state now the main result of this section:

\begin{theorem}
\label{thmmtr} The consistency strength of the statement $\mathfrak{r}_\kappa<2^\kappa=\kappa^{++}$ where $\kappa$ is weakly compact is at most the existence of a measurable cardinal $\kappa$ such that $o(\kappa)=\kappa^{+3}$.
\end{theorem}

\par\noindent\emph{Proof}. \newline
Let $\kappa$ be measurable with $o(\kappa)=\kappa^{+3}$, and assume that $2^\kappa=\kappa^{++}$.
Let $\bar{E}=(E_\eta:\eta\in\kappa^{+3})$ be a $\lhd$-increasing sequence of extenders over $\kappa$ and let $\mathbb{P}_{\bar{E}}$ be the associated extender-based Radin forcing.
Fix a sufficiently large regular cardinal $\chi>\kappa^{+3}$ and choose an increasing and continuous sequence $(N_i:i\leq\kappa^+)$ of $\bar{E}$-nice models.
Let $\delta_i=N_i\cap\kappa^{+3}$ be the characteristic ordinal of $N_i$ for every $i\leq\kappa^+$ and denote $\delta_{\kappa^+}$ by $\delta$.
Set $\hat{E}=\bar{E}\upharpoonright\delta=(E_\eta:\eta\in\delta)$.
We shall force $\mathfrak{r}_\kappa<2^\kappa=\kappa^{++}$ with $\mathbb{P}_{\hat{E}}$.
Fix a $V$-generic set $G\subseteq\mathbb{P}_{\hat{E}}$.
For every $i\in\kappa^+$ let $G_i=G\cap N_i$.
By Lemma \ref{lemsubforcing}, $G_i$ is equivalent to a $\mathbb{P}_{{\hat{E}}\upharpoonright\delta_i}$-generic set over $V$.

Recall that $E_{\delta_i}(\kappa)= \{A\subseteq\kappa:\kappa\in\jmath_{E_{\delta_i}}(A)\}$ is a normal measure over $\kappa$ for every $i\in\kappa^+$.
Working in $V[G_i]$ we define $\mathscr{U}_i$ as follows.
If $A$ is a subset of $\kappa$ in $V[G_i]$ then $A\in\mathscr{U}_i$ iff there exists a condition $p\in G_i$ such that $p_\rightarrow=(f,T),d={\rm dom}f$ and the following holds:
$$
\jmath_{E_{\delta_i}}(p)_{\langle{\rm mc}_{\delta_i}(d)\rangle}
\Vdash\check{\kappa}\in\jmath_{E_{\delta_i}}(\name{A}).
$$
One can show that $\mathscr{U}_i$ is a normal measure over $\kappa$ in $V[G_i]$, and it extends $E_{\delta_i}(\kappa)$.
A full detailed proof appears in \cite[Proposition 2.7]{MR3436372}.

Recall that $G^\alpha=\bigcup\{f^{p_\rightarrow}(\bar{\alpha}):p\in G,\bar{\alpha}\in{\rm dom}(f^{p_\rightarrow})\}$ is defined in $V[G]$, and the associated sequence is defined by $C^\alpha=\{\bar{\nu}_0:\bar{\nu}\in G^\alpha\}$.
We are interested mainly in the Radin club $C^\kappa$.
For every $\tau\in C^\kappa$ there is a condition $p\in G$ which forced it into $C^\kappa$.
This means that for some $\nu$ in the first level of the tree $T^{p_\rightarrow}$ we have $\tau=\nu(\bar{\kappa})_0$ and $p_{\langle\nu\rangle}\in G$.
We define $o^G(\tau)=o(\nu(\bar{\kappa}))$ whenever $\tau,p,\nu$ satisfy the above requirements, where $o(\nu(\bar{\kappa}))$ is the length of the sequence of extenders which appears in $\nu(\bar{\kappa})$.
Likewise, if $\tau\in C^\kappa$ and $\bar{\alpha}\in{\rm dom}(\nu)$ for some $p,\nu$ as above then we define $t_\alpha(\tau)=\nu(\bar{\alpha})_0$.
Now for every $i\in\kappa^+$ we define the following set:
$$
b_i = \{\mu\in C^\kappa: \mu\in{\rm dom}(t_{\delta_i}), o^G(\mu)=t_{\delta_i}(\mu)\}.
$$
The crucial point, proved in \cite[Proposition 2.12]{MR3436372}, is that $b_i$ generates the ultrafilter $\mathscr{U}_i$.
That is, if $A\in V[G_i]\cap\mathscr{U}_i$ then $V[G]\models b_i\subseteq^* A$.
For this, hark back to the definition of $\mathscr{U}_i$ which guarantees that there exists a condition $p\in G_i$ such that for every $\langle\nu\rangle\in T^{p_\rightarrow}$, if $p_{\langle\nu\rangle}\Vdash t_{\delta_i}(\nu(\bar{\kappa})_0)=o^{\name{G}}(\nu(\bar{\kappa})_0)$ then $p_{\langle\nu\rangle}\Vdash\nu(\bar{\kappa})_0\in\name{A}$.

Working in $V[G]$, let $\mathscr{R}=\{b_i:i\in\kappa^+\}$.
We claim that $\mathscr{R}$ is an $\mathfrak{r}_\kappa$-family, thus proving that $\kappa^+ = \mathfrak{r}_\kappa<2^\kappa=\kappa^{++}$ in $V[G]$.
To see this, suppose that $B\subseteq\kappa, B\in V[G]$.
From Lemma \ref{lemproper} we know that $B\in V[G_i]$ for some $i\in\kappa^+$.
Since $\mathscr{U}_i$ is an ultrafilter over $\kappa$ in $V[G_i]$ we see that $B\in\mathscr{U}_i$ or $(\kappa-B)\in\mathscr{U}_i$.
By the crucial point mentioned above, either $b_i\subseteq^* B$ or $b_i\subseteq^*(\kappa-B)$, so $B$ fails to split all the elements of $\mathscr{R}$.
But $B$ was an arbitrary subset of $\kappa$ in $V[G]$, so $\mathscr{R}$ is an $\mathfrak{r}_\kappa$-family.

\hfill \qedref{thmmtr}

In the above model we see that $\mathfrak{s}_\kappa>\kappa$ holds in $V[G]$,
so $\kappa$ is weakly compact in the generic extension.
We do not know how to obtain a similar statement at strongly inaccessible cardinals which are not weakly compact.
Likewise, we do not know how to get $\mathfrak{r}_\kappa<2^\kappa=\kappa^{++}$ at successor cardinals, see \cite{RagSh1}.

Remark that one can avoid the assumption $2^\kappa=\kappa^{++}$ in the ground model, and obtain a similar result.
Indeed, by rendering the same forcing construction, we necessarily increase $2^\kappa$ to $\kappa^{++}$ while creating local repeat points even if we assume $2^\kappa=\kappa^+$ in the ground model.
However, the above construction seems more transparent.

We do not know how to obtain a lower bound with respect to the consistency strength of $\mathfrak{r}_\kappa<2^\kappa$.
The analogy between $\mathfrak{r}_\kappa$ and $\mathfrak{s}_\kappa$ is suggestive, and one tends to think that large cardinals are required.
On the other hand, \cite{RagSh1} shows that a measurable cardinal is sufficient for $\mathfrak{r}_\kappa<2^\kappa$ or even $\mathfrak{u}_\kappa<2^\kappa$.

\begin{question}
\label{qr} What is the exact consistency strength of the statement $\mathfrak{r}_\kappa<2^\kappa$? How about the statement $\mathfrak{r}_\kappa<2^\kappa=\kappa^{++}$?
\end{question}

\newpage

\section{Measurable but not supercompact}

In this section we deal with several questions which apparently require more than a measurable cardinal $\kappa$ with $o(\kappa)=\kappa^{+3}$.
Basically, we wish to prove that $\mathfrak{r}_\kappa<2^\kappa=\kappa^{++}$ or $\mathfrak{s}_\kappa>\kappa^+$ while preserving the measurability of $\kappa$.
This can be done if $\kappa$ is supercompact.
However, we would like to get this inequality from hyper-measurability.

We shall address the following questions:
\begin{itemize}
\item What is the consistency strength of the statement $\mathfrak{r}_\kappa<2^\kappa=\kappa^{++}$ where $\kappa$ is a measurable cardinal (see also \cite{RagSh1}, Question 17)?
\item What is the consistency strength of the statement $\mathfrak{s}_\kappa=\kappa^{++}$ where $\kappa$ is a measurable cardinal (\cite{MR3436372}, Question 3.2)?
\item Is it consistent that $\kappa$ is measurable but not supercompact, and $\binom{2^\kappa}{\kappa}\rightarrow\binom{2^\kappa}{\kappa}$ holds (\cite{MR3509813}, Question 4.3)?
\end{itemize}

A few words about the ultrafilter number are in order.
Recall that $\mathfrak{u}_\kappa$ is the minimal size of a base $\mathcal{B}$ of some uniform ultrafilter $\mathscr{U}$ over $\kappa$, where a base $\mathcal{B}$ of $\mathscr{U}$ satisfies $\mathcal{B}\subseteq\mathscr{U}$ and for every $A\in\mathscr{U}$ there is $B\in\mathcal{B}$ such that $B\subseteq A$.
The reaping number $\mathfrak{r}_\kappa$ is very close to $\mathfrak{u}_\kappa$.
Easily $\mathfrak{r}_\kappa\leq\mathfrak{u}_\kappa$ and the question whether $\mathfrak{r}_\kappa<\mathfrak{u}_\kappa$ is consistent for some $\kappa>\aleph_0$ is widely open.

We have seen that $o(\kappa)=\kappa^{+3}$ is sufficient for the consistency of $\mathfrak{r}_\kappa<2^\kappa=\kappa^{++}$, and
we do not know what is the consistency strength of $\mathfrak{u}_\kappa<2^\kappa=\kappa^{++}$.
This statement has been proved by Gitik and Shelah in \cite{MR1632081} by starting from a huge cardinal in the ground model.
Hugeness has been replaced by supercompactness in \cite{MR3564375} and \cite{MR3201820}.
Unfortunately, it seems that the methods of the previous section are not applicable to $\mathfrak{u}_\kappa$.

Let us begin with some preliminaries.
Suppose that $\kappa$ is measurable and $o(\kappa)=\kappa^{+3}$.
Let $\bar{E}=(E_\eta:\eta\in\kappa^{+3})$ be a $\lhd$-increasing sequence of extenders over $\kappa$.
Recall that if $d\in[\kappa^{+3}]^{\leq\kappa}$ then $E(d) = \bigcap\{E_\eta(d):\eta\in\kappa^{+3}\}$.
An ordinal $\zeta\in\kappa^{+3}$ is called \emph{a repeat point} of $\bar{E}$ iff $\bigcap\{E_\eta(d):\eta<\zeta\}=E(d)$ for every $d\in[\kappa^{+3}]^{\leq\kappa}$.
The following definition comes from \cite{MR3436372}, and will be central in our arguments:

\begin{definition}
\label{deflocalrp} Local repeat points. \newline
Let $\bar{E} = (E_\eta:\eta\in\lambda)$ be a $\lhd$-increasing sequence of extenders over $\kappa$. \newline
An ordinal $\delta\in\lambda$ is a local repeat point of $\bar{E}$ iff $\bigcap\{E_\eta(d):\eta<\delta\}=E(d)$ for every $d\in[\delta]^{\leq\kappa}$.
\end{definition}

The locality is reflected in the fact that the repeatedness is required only for subsets of $\delta$.
We introduce now a general method for creating local repeat points.
Assume that $\chi>\kappa^{+3}$ is a sufficiently large regular cardinal.
We shall say that $N$ is $\bar{E}$-nice iff $N\prec\mathcal{H}(\chi), \kappa^{++}\subseteq N, {}^\kappa N\subseteq N, |N|=\kappa^{++}$ and $\bar{E}\in N$.
Notice that if $N$ is $\bar{E}$-nice then $\delta_N=N\cap\kappa^{+3}$ is an ordinal in $\kappa^{+3}$, called \emph{the characteristic ordinal of $N$}.

\begin{lemma}
\label{lemlocalrp} Assume that $2^\kappa=\kappa^{++}$. \newline
If $N$ is $\bar{E}$-nice then $\delta_N$ is a local repeat point of $\bar{E}$.
\end{lemma}

\par\noindent\emph{Proof}. \newline
Fix an ordinal $\beta\in\kappa^{+3}$ and a set $d\in[\kappa^{+3}]^{\leq\kappa}$.
Consider the intersection $F_\beta(d)=\bigcap\{E_\eta(d):\eta<\beta\}$.
If $\beta<\gamma<\kappa^{+3}$ then $F_\beta(d)\supseteq F_\gamma(d)$, hence $(F_\beta(d):\beta\in\kappa^{+3})$ forms a $\subseteq$-decreasing sequence of filters over $V_\kappa$.
Since $|\mathcal{P}(V_\kappa)|=2^\kappa=\kappa^{++}$, there must be an ordinal $\gamma=\gamma(d)$ for which $\bigcap\{E_\eta(d):\eta<\gamma\}=E(d)$.

We focus now on $\delta_N=N\cap\kappa^{+3}$.
Since ${}^\kappa N\subseteq N$, if $d\in[\delta_N]^{\leq\kappa}$ then $d\in N$.
For such $d$ there is an ordinal $\gamma(d)$ as described above, and by elementarity we may choose $\gamma(d)\in N$, so $\gamma(d)\in N\cap\kappa^{+3}=\delta_N$ whenever $d\in[\delta_N]^{\leq\kappa}$.
It follows that $\bigcap\{E_\eta(d):\eta<\delta_N\}=E(d)$ for every $d\in[\delta_N]^{\leq\kappa}$, so $\delta_N$ is a local repeat point of $\bar{E}$, as required.

\hfill \qedref{lemlocalrp}

Let $N$ be $\bar{E}$-nice and consider $\mathbb{P}_{\bar{E}}\cap N$.
Denote the sequence $(E_\eta:\eta<\delta_N)$ by $\bar{E}\upharpoonright\delta_N$.
The following lemma says that forcing with $\mathbb{P}_{\bar{E}}\cap N$ is essentially like forcing with $\mathbb{P}_{{\bar{E}}\upharpoonright\delta_N}$.
The two parts of the lemma are proved in \cite[Lemma 2.1]{MR3436372} and \cite[Proposition 2.4]{MR3436372} respectively.

\begin{lemma}
\label{lemsubforcing} Suppose that $N$ is $\bar{E}$-nice.
\begin{enumerate}
\item [$(\aleph)$] $N\cap\mathbb{P}_{\bar{E}}$ is a complete subforcing of $\mathbb{P}_{\bar{E}}$.
\item [$(\beth)$] $N\cap\mathbb{P}_{\bar{E}}$ is isomorphic to $\mathbb{P}_{{\bar{E}}\upharpoonright\delta_N}$.
\end{enumerate}
\end{lemma}

\hfill \qedref{lemsubforcing}

Recall that if $o(\kappa)=\lambda$ then there is a $\lhd$-increasing sequence of extenders $\bar{E}=(E_\alpha:\alpha\in\lambda)$ such that $\kappa={\rm crit}(\bar{E})$, the sequence of ordinals $(\sigma(E_\alpha):\alpha\in\lambda)$ is strictly increasing and $\sigma(E_\alpha)<\lambda$ for every $\alpha\in\lambda$.
Given such a sequence $\bar{E}$ we define $C_{\bar{E}} = \{\alpha\in\lambda:\forall\beta<\alpha,\sigma(E_\beta)<\alpha\}$.
We are assuming that $\lambda$ is a regular cardinal and hence $C_{\bar{E}}$ is a club subset of $\lambda$.

In the following couple of lemmata we prove that one can refine the assumptions on the sequence $\bar{E}$.
The first lemma requires $\lambda$ to be a regular cardinal, and the second lemma employs the weak compactness of $\lambda$.

\begin{lemma}
\label{lemheight} Suppose that $\kappa<\lambda, \kappa$ is a measurable cardinal, $o(\kappa)=\lambda$ and $\lambda$ is regular. \newline
Then there is a $\lhd$-increasing sequence of extenders $\bar{E}=(E_\alpha:\alpha\in\lambda)$ such that $\kappa={\rm crit}(\bar{E})$ and $\sigma(E_\delta)=\delta$ whenever $\delta\in C_{\bar{E}}$ is strongly inaccessible.
\end{lemma}

\par\noindent\emph{Proof}. \newline
Begin with an extender-sequence $\bar{e}=(e_\alpha:\alpha\in\lambda)$ whose critical point is $\kappa$.
We always assume that $\delta\leq\sigma(e_\delta)$ but it may happen that $\delta<\sigma(e_\delta)$.
In order to obtain $\delta=\sigma(e_\delta)$ at many points we shall replace each $e_\alpha$ by $E_\alpha$ using the following strategy.

Firstly, we choose a function $f:\kappa\rightarrow V_\kappa$ such that $\bar{e}\upharpoonright\alpha=\jmath_{e_\alpha}(f)(\kappa)$ for every $\alpha\in\lambda$.
Secondly, let $\delta\in C_{\bar{e}}$ be strongly inaccessible.
Let $e_\delta\upharpoonright\delta$ be \emph{the cutback} of $e_\delta$, namely the extender obtained by taking only the measures of $e_\delta$ which correspond to generators below $\delta$.
Notice that $\sigma(e_\delta\upharpoonright\delta)=\delta$.
Finally, define $E_\delta=e_\delta\upharpoonright\delta$ for every strongly inaccessible $\delta\in C_{\bar{e}}$ and let $E_\delta=e_\delta$ otherwise.
We claim that the sequence $\bar{E}=(E_\alpha:\alpha\in\lambda)$ satisfies our lemma.

We must check that $\bar{E}$ is still $\lhd$-increasing.
Let $\jmath_{e_\delta}:V\rightarrow M_{e_\delta}$ and $\jmath_{E_\delta}:V\rightarrow M_{E_\delta}$ be the canonical ultrapower embeddings.
We define the natural mapping $k_\delta:M_{E_\delta}\rightarrow M_{e_\delta}$ by $k_\delta(\jmath_{E_\delta}(f)(\alpha))=\jmath_{e_\delta}(f)(\alpha)$ whenever $f\in V, \kappa={\rm dom}f$ and $\alpha\in\delta$.
This definition yields the following commuting diagram:

$$
\xymatrix{
{V} \ar[ddrr]_{\jmath_{e_\delta}} \ar[rr]^{\jmath_{E_\delta}} & & M_{E_\delta} \ar[dd]^{k_\delta} \\ \\
& & M_{e_\delta} }
$$

Our argument will be based on the fact that ${\rm crit}(k_\delta)\geq\delta$.
To verify this, let $f:\kappa\rightarrow\kappa$ be the identity function, so $f\in V$ and $k_\delta(\alpha)=k_\delta(\jmath_{E_\delta}(f)(\alpha)) = (\jmath_{e_\delta}(f)(\alpha))=\alpha$ for every $\alpha\in\delta$ by the definition of $k_\delta$ and the elementarity of $\jmath_{E_\delta}$ and $\jmath_{e_\delta}$.

Observe also that $k_\delta(\bar{E})=\bar{e}\upharpoonright\delta$.
Indeed, let $g:\kappa\rightarrow V_\kappa$ be defined by the requirement $\bar{e}\upharpoonright\delta=\jmath_{e_\delta}(g)(\kappa)$.
Letting $\bar{E}=\jmath_{E_\delta}(g)(\kappa)$ we see that $\bar{E}\in M_{E_\delta}$.
Apply $k_\delta$ to the equality $\bar{E}=\jmath_{E_\delta}(g)(\kappa)$, and infer that $k_\delta(\bar{E}) = \kappa_\delta(\jmath_{E_\delta}(g)(\kappa)) = \jmath_{e_\delta}(g)(\kappa) = \bar{e}\upharpoonright\delta$.
But $\delta\in C_{\bar{e}}$ and $\delta$ is strongly inaccessible, so $\bar{e}\upharpoonright\delta\subseteq V_\delta$.
Since $k_\delta(\bar{E})=\bar{e}\upharpoonright\delta\subseteq V_\delta$ and ${\rm crit}(k_\delta)\geq\delta$ we conclude that $\bar{E}=\bar{e}\upharpoonright\delta$.
In particular, $\bar{E}$ is $\lhd$-increasing, so we are done.

\hfill \qedref{lemheight}

The next lemma provides a stationary subset $S\subseteq\lambda$ (associated with $\bar{E}$) with the following property.
If $\delta,\varepsilon\in S, \delta<\varepsilon$ and both are local repeat points then the normal measures $\mathscr{U}_\delta$ and $\mathscr{U}_\varepsilon$ (described in the previous section) satisfy $\mathscr{U}_\delta\subseteq\mathscr{U}_\varepsilon$.
This essential property enables us to keep the measurability of $\kappa$ in the generic extension.
It will be useful when trying to force $\mathfrak{u}_\kappa<2^\kappa$, and moreover when trying to obtain a normal ultrafilter $\mathscr{U}$ over $\kappa$ such that ${\rm Ch}(\mathscr{U})=\kappa^+<2^\kappa$.
Likewise, it will be useful while forcing $\mathfrak{s}_\kappa>\kappa^+$ where $\kappa$ is measurable.
Both statements will be derived from much less than supercompactness.

Let $\bar{E}=(E_\alpha:\alpha\in\lambda)$ be as promised in Lemma \ref{lemheight}.
Let $\mathbb{P}_{\bar{E}}$ be the associated extender-based Radin forcing.
Fix a local repeat point $\delta<\lambda$.
Recall that $E_\delta(\kappa)=\{A\subseteq\kappa:\kappa\in\jmath_{E_\delta}(A)\}$ is normal.
Moreover, $\kappa$ remains measurable in the generic extension $V[G_{\bar{E}\upharpoonright\delta}]$.

To see this, we define a normal measure $\mathscr{U}_\delta$ over $\kappa$, in $V[G_{\bar{E}\upharpoonright\delta}]$, as follows.
If $\name{A}$ is a $V$-name for a subset of $\kappa$ in $V[G_{\bar{E}\upharpoonright\delta}]$ then $A\in\mathscr{U}_\delta$ iff there is a condition $p\in G_{\bar{E}\upharpoonright\delta}$ such that  $p_\rightarrow=(f,T), d={\rm dom}f$ and $\jmath_{E_\delta}(p)_{\langle{\rm mc}_\delta(d)\rangle}\Vdash \check{\kappa}\in \jmath_{\bar{E}\upharpoonright\delta}(\name{A})$.

\begin{lemma}
\label{lemweakcomp} Suppose that $\kappa$ is measurable with $o(\kappa)=\lambda$, and $\lambda$ is a weakly compact cardinal greater than $\kappa$. Let $\bar{E}$ and $C_{\bar{E}}$ be as in the previous lemma. \newline
Then there exists a set $S_{\bar{E}}\subseteq\lambda$ such that:
\begin{enumerate}
\item [$(\aleph)$] $S_{\bar{E}}$ is a stationary subset of $\lambda$.
\item [$(\beth)$] If $\delta_0,\delta_1\in S_{\bar{E}}, \delta_0<\delta_1$ and both are local repeat points then for every $V$-generic subset $G_{\bar{E}\upharpoonright\delta_1}\subseteq \mathbb{P}_{\bar{E}\upharpoonright\delta_1}$ we have $\mathscr{U}_{\delta_0}=\mathscr{U}_{\delta_1}\cap V[G_{\bar{E}\upharpoonright\delta_0}]$, where $G_{\bar{E}\upharpoonright\delta_0}$ is the restriction of $G_{\bar{E}\upharpoonright\delta_1}$ to $\delta_0$.
\end{enumerate}
\end{lemma}

\par\noindent\emph{Proof}. \newline
Let $\chi$ be a sufficiently large regular cardinal above $\lambda$.
Let $M\prec(\mathcal{H},\in\bar{E})$ be a model of size $\lambda$ such that ${}^{<\lambda}M\subseteq M$.
Since $\lambda$ is weakly compact, there exists a $\lambda$-complete normal filter $\mathcal{F}$ over $\lambda$ such that if $A\in M\cap\mathcal{P}(\lambda)$ then $A\in\mathcal{F}\vee(\lambda-A)\in\mathcal{F}$.

Let $N$ be the ultrapower formed by $\mathcal{F}$ with respect to $M$, and let $\pi:M\rightarrow N$ be the canonical embedding.
Let $\bar{E}^N=(E^N_\alpha:\alpha\in\pi(\lambda))$ be the extender sequence $\pi(\bar{E})$.
By elementarity, $\bar{E}^N$ is $\lhd$-increasing.
Since $\lambda={\rm crit}(\pi)$ we see that $\bar{E}^N\upharpoonright\lambda=\bar{E}$.
If $\delta\in C_{\bar{E}}$ is strongly inaccessible then we are assuming that $\sigma(E_\delta)=\delta$ (by the previous lemma) and hence $\sigma(E^N_\lambda)=\lambda$ using the fact that $\mathcal{F}$ is normal.
Define:

$$
S_{\bar{E}} = \{\delta\in C_{\bar{E}}:\delta\ \text{is inaccessible and}\
\forall d\in[\delta]^{\leq\kappa}, E_\delta(d)=E^N_\lambda(d)\}.
$$
We claim that $S_{\bar{E}}$ satisfies the statements of the lemma.

The fact that $S_{\bar{E}}$ is stationary follows, again, from the normality of $\mathcal{F}$.
Suppose that $\delta_0,\delta_1\in S_{\bar{E}}$ and $\delta_0<\delta_1$.
If $d\in[\delta_0]^{\leq\kappa}$ then $d\in[\delta_1]^{\leq\kappa}$ as well, and $E_{\delta_0}(d)=E_{\delta_1}(d)$ since both equal to $E^N_\lambda(d)$.
Hence if both $\delta_0$ and $\delta_1$ are local repeat points then $\mathscr{U}_{\delta_0}\subseteq\mathscr{U}_{\delta_1}$.
This means that $\mathscr{U}_{\delta_0}\subseteq\mathscr{U}_{\delta_1}\cap V[G_{\bar{E}\upharpoonright\delta_0}]$, and since $\mathscr{U}_{\delta_0}$ is an ultrafilter in $V[G_{\bar{E}\upharpoonright\delta_0}]$ we see that $\mathscr{U}_{\delta_0}=\mathscr{U}_{\delta_1}\cap V[G_{\bar{E}\upharpoonright\delta_0}]$, as desired.

\hfill \qedref{lemweakcomp}

We can prove now the consistency of $\mathfrak{r}_\kappa<2^\kappa=\kappa^{++}$ where $\kappa$ is measurable but not supercompact.

\begin{theorem}
\label{thm1ucom} It is consistent that $\kappa$ is measurable and $\mathfrak{r}_\kappa<2^\kappa$ relative to a measurable cardinal $\kappa$ such that $o(\kappa)=\lambda$ and $\lambda$ is a weakly compact cardinal greater than $\kappa$.
\end{theorem}

\par\noindent\emph{Proof}. \newline
Let $\kappa$ be a measurable cardinal, $\lambda>\kappa$, and assume that $\lambda$ is weakly compact and $o(\kappa)=\lambda$ in the ground model.
We assume that the GCH holds in the ground model.
Let $\bar{E}=(E_\alpha:\alpha\in\lambda)$ be a $\lhd$-increasing sequence of extenders which satisfies the above lemmata.
Namely, $\kappa={\rm crit}(\bar{E})$, if $\delta\in C_{\bar{E}}$ is strongly inaccessible then $\sigma(E_\delta)=\delta$ and $S_{\bar{E}}$ is a stationary subset of $\lambda$ in which local repeat points give rise to an increasing sequence of measures with respect to inclusion.

Let $\chi$ be a sufficiently large regular cardinal above $\lambda$.
An elementary submodel $N\prec\mathcal{H}(\chi)$ will be called $\bar{E}$-suitable iff $N\cap\lambda=\delta_N>\kappa^{++}, |N|=|\delta_N|<\lambda, {}^\kappa N\subseteq N$ and $\bar{E}\in N$.
This definition resembles the definition of $\bar{E}$-nice models, but here the characteristic ordinal $\delta_N$ and the cardinality of $N$ are strictly above $\kappa^{++}$.

Consider an increasing continuous sequence $(N_i:i\in\lambda)$ of $\bar{E}$-suitable elementary submodels of $\mathcal{H}(\chi)$.
The corresponding sequence of ordinals $(\delta_{N_i}:i\in\lambda)$ forms a club subset of $\lambda$.
Hence $S_{\bar{E}}\cap\{\delta_{N_i}:i\in\lambda\}$ is a stationary subset of $\lambda$.
We can choose, therefore, a shorter increasing continuous sequence $(N_i:i\in\kappa^+)$ of $\bar{E}$-suitable structures such that $\delta_{N_{i+1}}\in S_{\bar{E}}$ for every $i\in\kappa^+$.

Denote $\delta_{N_i}$ by $\delta_i$ for every $i\in\kappa^+$ and let $\delta=\delta_N$ where $N=\bigcup_{i\in\kappa^+}N_i$.
Let $\hat{E}=\bar{E}\upharpoonright\delta$ and $\mathbb{P}_{\hat{E}}=\mathbb{P}_{\bar{E}}\cap N$.
We shall force over $V$ with $\mathbb{P}_{\hat{E}}$, and we claim that if $G\subseteq\mathbb{P}_{\hat{E}}$ is generic then $V[G]$ satisfies the following statements:
\begin{enumerate}
\item [$(a)$] $2^\kappa=\delta>\kappa^+$.
\item [$(b)$] $\mathfrak{r}_{\kappa}=\kappa^+$.
\end{enumerate}
For part $(a)$ recall that $G^\alpha=\bigcup\{f^{p_\rightarrow}(\bar{\alpha}): p\in G, \bar{\alpha}\in{\rm dom}(f^{p_\rightarrow})\}$ and $C^\alpha=\{\bar{\nu}_0:\bar{\nu}\in G^\alpha\}$ for every $\alpha\in[\kappa,\delta)$.
Since $\alpha\neq\beta\Rightarrow C^\alpha\neq C^\beta$ one concludes that $V[G]\models 2^\kappa\geq\delta$.
Actually, $V[G]\models 2^\kappa=\delta$ since $\mathbb{P}_{\hat{E}}$ is $\delta$-cc.
Of course, $\delta>\kappa^+$ by its definition.

For part $(b)$ we indicate that Lemma \ref{lemsubforcing} holds true when applied to $\bar{E}$-suitable models (the proof in \cite{MR3436372} is phrased with respect to $\bar{E}$-nice models, but the distinction between the concepts has no influence on the proof).
In particular, $N_{i+1}\cap\mathbb{P}_{\hat{E}}$ is a complete subforcing of $\mathbb{P}_{\hat{E}}$ for every $i\in\kappa^+$.
Aiming to show that $\mathfrak{u}_\kappa^{\rm com}=\kappa^+$ we define a $\kappa$-complete ultrafilter $\mathscr{U}$ over $\kappa$ in $V[G]$ as follows:
$$
\mathscr{U}=\bigcup\{\mathscr{U}_{\delta_{i+1}}:i\in\kappa^+\}.
$$
If $i\in\kappa^+$ then $\delta_{i+1}\in S_{\bar{E}}$ and hence it is a local repeat point of $\bar{E}$.
Since $\hat{E}=\bar{E}\upharpoonright\delta$ we see that $\delta_{i+1}$ is a local repeat point of $\hat{E}$ as well.
It follows that $\mathscr{U}_{\delta_{i+1}}$ is a $\kappa$-complete (and normal) ultrafilter over $\kappa$ in the intermediate extension $V[G\cap N_{i+1}]$.

If $A\in V[G]\cap\mathcal{P}(\kappa)$ then $A\in V[G\cap N_j]$ for some $j\in\kappa^+$ and hence $A\in V[G\cap N_i]$ for every $i\in[j,\kappa^+)$.
Without loss of generality $A\in\mathscr{U}_{\delta_j}$, so $A\in\mathscr{U}_{\delta_i}$ for every $i\in[j,\kappa^+)$ and hence $A\in\mathscr{U}$.
This argument shows that $\mathscr{U}$ is an ultrafilter over $\kappa$.
A similar argument shows that $\mathscr{U}$ is $\kappa$-complete.
In particular, $\kappa$ remains measurable in the generic extension.
Let $\mathscr{B}=\{b_i:i\in\kappa^+\}$.
Notice that $\mathscr{B}$ is unsplittable and hence $\mathfrak{r}_\kappa=\kappa^+<2^\kappa=\delta$, as required.

\hfill \qedref{thm1ucom}

In the above theorem we obtain $\kappa^+=\mathfrak{r}_\kappa<2^\kappa=\delta$, where $\delta$ is quite large.
However, one can get a similar statement with $2^\kappa=\kappa^{++}$.
For this end, force as in the above theorem and then collapse $\delta$ to $\kappa^{++}$.
The completeness of the collapse implies that no new subsets of $\kappa$ are introduced, and our result follows.
This idea seems inapplicable to the model of \cite{RagSh1}.

We also mention here an open problem from \cite{MR3509813}, for which we can give a positive answer.
We know that one can force $\binom{2^\kappa}{\kappa}\rightarrow\binom{2^\kappa}{\kappa}_2$ at a supercompact cardinal $\kappa$.
Question 4.3 in \cite{MR3509813} is whether one can force this relation at a measurable but not supercompact cardinal.

\begin{corollary}
\label{cor1012} It is consistent that the positive relation $\binom{2^\kappa}{\kappa}\rightarrow\binom{2^\kappa}{\kappa}_2$ holds at a measurable cardinal $\kappa$ which is not supercompact.
\end{corollary}

\par\noindent\emph{Proof}. \newline
Assume that there is a measurable cardinal $\kappa$ such that $o(\kappa)=\lambda$ and $\lambda>\kappa$ is weakly compact.
Assume, further, that there are no supercompact cardinals in the ground model.
Use the above theorem to force $\mathfrak{r}_\kappa<2^\kappa$.
In the generic extension we have $\binom{2^\kappa}{\kappa}\rightarrow\binom{2^\kappa}{\kappa}_2$, as required.

\hfill \qedref{cor1012}

It seems, however, that our method cannot be applied to strongly inaccessible but not weakly compact cardinals, see Question 4.2 of \cite{MR3509813} and the discussion at the end of this section.

We conclude this section with the splitting number.
Question 3.2 of \cite{MR3436372} is about the consistency strength of $\mathfrak{s}_\kappa>\kappa^+$ where $\kappa$ is measurable.
Similar methods to those employed in the theorem above will provide a proof to the following:

\begin{theorem}
\label{thm2skappa} It is consistent that $\kappa$ is measurable and $\mathfrak{s}_\kappa>\kappa^+$ relative to a measurable cardinal $\kappa$ such that $o(\kappa)=\lambda$ and $\lambda$ is a weakly compact cardinal greater than $\kappa$.
\end{theorem}

\par\noindent\emph{Proof}. \newline
Let $\kappa$ be measurable in $V$, and let $\lambda$ be a weakly compact cardinal strictly above $\kappa$ so that $o(\kappa)=\lambda$.
We shall assume that the GCH holds in $V$.
Let $\bar{E},C_{\bar{E}},S_{\bar{E}}$ be as in the proof of Theorem \ref{thm1ucom}.

Choose an increasing continuous sequence $(N_i:i\in\kappa^{++})$ of $\bar{E}$-suitable models with $\delta_{i+1}\in S_{\bar{E}}$ for every $i\in\kappa^{++}$.
Let $\hat{E}$ be $\bar{E}\upharpoonright\delta_{\kappa^{++}}$.
As before, we shall force with $\mathbb{P}_{\hat{E}}$ over $V$, so we choose a $V$-generic subset $G\subseteq\mathbb{P}_{\hat{E}}$.
Let $\mathscr{U}_{\delta_{i+1}}$ be the normal measure over $\kappa$ in $V[G\cap N_{i+1}]$, and let $(b_i:i\in\kappa^{++})$ be the corresponding sequence of subsets of $\kappa$ defined within the proof of Theorem \ref{thm1ucom}.

Suppose that $\mathcal{F}\subseteq[\kappa]^\kappa\cap V[G]$ and $|\mathcal{F}|\leq\kappa^+$.
Toward contradiction assume that $\mathcal{F}$ is a splitting family.
Enumerate the elements of $\mathcal{F}$ by $\{x_\alpha:\alpha\in\kappa^+\}$, using repetitions if needed.
For every $\alpha\in\kappa^+$ choose $j_\alpha\in\kappa^{++}$ so that if $i\in[j_\alpha,\kappa^{++})$ then $x_\alpha\in V[G\cap N_i]$.
Let $j=\bigcup_{\alpha\in\kappa^+}j_\alpha$, so $j\in\kappa^{++}$ and $\mathcal{F}\in V[G\cap N_{j+1}]$.

If $\mathcal{S}\subseteq[\kappa]^\kappa$ is a splitting family then every $x\in\mathcal{S}$ can be replaced by $\kappa-x$ keeping the splitting property of the family.
Hence if $\mathscr{U}$ is an ultrafilter over $\kappa$ then one can assume that $\mathcal{S}\subseteq\mathscr{U}$.
In particular, we may assume that $\mathcal{F}\subseteq\mathscr{U}_{\delta_{j+1}}$.
Consequently, $b_{\delta_{j+1}}\subseteq^*x_\alpha$ for every $\alpha\in\kappa^+$ and hence $\mathcal{F}$ is not a splitting family, a contradiction.

We may conclude from the previous paragraph that $\mathfrak{s}_\kappa>\kappa^+$.
Likewise, we know that $\kappa$ is measurable in $V[G]$.
To see this, define $\mathscr{U}=\bigcup\{\mathscr{U}_{\delta_{i+1}}:i\in\kappa^{++}\}$.
In Theorem \ref{thm1ucom} we have seen that $\mathscr{U}$ is a $\kappa$-complete ultrafilter over $\kappa$ in $V[G]$, so $\kappa$ is measurable in $V[G]$ and our proof is accomplished.

\hfill \qedref{thm2skappa}

The methods of this section for proving the consistency of $\mathfrak{u}_\kappa<2^\kappa$ are based on increasing $\mathfrak{s}_\kappa$.
If $\kappa$ is not weakly compact then $\mathfrak{s}_\kappa\leq\kappa$, hence out of the scope of the above approach.
The method of \cite{RagSh1} is applicable to accessible cardinals.
This leaves open the case of strongly inaccessible but not weakly compact cardinals:

\begin{question}
\label{qinacc} Is it consistent that $\kappa$ is strongly inaccessible but not weakly compact, and $\mathfrak{r}_\kappa<2^\kappa$?
\end{question}

\newpage

\bibliographystyle{amsplain}
\bibliography{arlist}

\end{document}